\documentclass[12pt]{elsarticle}
\begin{document}
%
%
\noindent
{\Large \bf THE WAVE EQUATION FOR LEVEL SETS IS
\\[0.5cm] NOT A HUYGENS' EQUATION}
\\[1cm]

{\bf WOLFGANG QUAPP} 
\\
{Mathematisches Institut, Universit\"{a}t Leipzig, PF 100920, D-04009 Leipzig,\\ Germany, \ \
quapp@uni-leipzig.de,\ \ corresponding author, tel.: 49-(0)341-97-32162}
\\

{\bf JOSEP MARIA BOFILL}
\\
{Departament de Qu\'{i}mica Org\`{a}nica, Universitat de Barcelona; 
Institut de Qu\'{i}mica Te\`{o}rica i Computacional,  Universitat de Barcelona, (IQTCUB), 
Mart\'{i} i Franqu\`{e}s, 1, 08028 Barcelona, Spain, \ \ jmbofill@ub.edu}
\\

September 12, 2013
\\
%

\noindent
{\bf Abstract:}\\
Any surface can be foliated into equipotential hypersurfaces of the level sets. 
A current result is that the contours are the progressing wave fronts of a
certain hyperbolic partial differential equation, a wave equation. 
It is connected with the gradient lines, as well as with a corresponding eikonal equation.
The level of a surface point, seen as an additional coordinate, plays the central role in this treatment. 
A wave solution can be a sharp front. 
Here the validity of the Huygens' principle (HP) is of interest: 
there is no wake of the wave solutions in every dimension, 
if a special Cauchy initial value problem is posed. 
Additionally, there is no distinction into odd or even dimensions.
To compare this with Hadamard's 'minor premise' for a strong HP,
we calculate differential geometric objects like Christoffel symbols, curvature tensors and  
geodesic lines, to test the validity of the strong HP. 
However, for the differential equation for level sets, the main criteria are not fulfilled for the strong HP in the 
sense of Hadamard's 'minor premise'. 
\\

\noindent
Keywords{: 
 Contours; \ steepest ascent; \ 
wave equation; \ progressing waves; 
\\   
\hspace*{2.0cm} 
Huygens' principle. }

\section{Introduction}
\label{Kap1}
In this paper we treat a hypersurface in an  $I\!\!R^{N+1}$. It will be described by a unique function  $\nu=V({\bf q})$
where ${\bf q}=(q^1,...,q^N)^T \in I\!\!R^{N}$ are space coordinates and $\nu$ is the (N+first) coordinate. 
Usually in applications  \cite{mez,hei,wales}, 
the region of interest of ${\bf q}$ is not the full $I\!\!R^{N}$ but is only a certain, local region.
The surface should be a continuous function with respect to the coordinates. It should also have
continuous derivatives up to a certain order not specified here, 
but required by the operations which are to be carried out.

The steepest descent lines are orthogonal trajectories to the contour hypersurfaces,
$V({\bf q})=\nu =constant$, if the corresponding metric relations are used \cite{quhe84,wqrev}.
In this paper we will assume $N$ orthogonal and equidistant coordinates, ${\bf q}=(q^1,...,q^N)^T$.
Then the metric matrix in $I\!\!R^N$ reduces to the unity matrix and we have a Euclidean behaviour.
\\

From a theoretical point of view (however not from a numerical one) the steepest descent lines
and the inverse ones, the steepest ascent lines are equivalent.
The steepest ascent/steepest descent lines emerging from a minimum or a saddle point of the surface V
can be seen as travelling in an orthogonal manner through the contour hypersurfaces
of this surface V.  It should be noted that the construction of the
contour  hypersurface, $V({\bf q})= \nu=constant$, is such that all points satisfying this
equation possess the same equipotential difference with respect to another
contour hypersurface.
It is similar to the construction of Fermat-Huygens of the propagation of wave fronts and rays. 
Note that this construction and the Hamilton-Jacobi theory are connected \cite{cour}.
Using the analogy, we have proposed a wave equation for contour hypersurfaces of the surface V \cite{boqu12},
which we report in Section \ref{Kap2}. Some preliminar work to the idea was reported in Ref.\,\cite{mog08}. 
The theory of wave equations in an $(N+1)$-dimensional 
space is complicated. 
One has to pose the Cauchy initial value problem (IVP), 
and the general theory \cite{fried,guen,guen91,hadam} leads to (local) integral formulas over regions of 
the so-called characteristic conoid. 
\\

``It is a familiar fact from daily life that propagation of waves is very different in 2 and 3 dimensions. 
When a pebble falls in water at a point $P$ the initial ripple on a circle around $P$ will be followed by subsequent ripples. 
Thus a given point $Q$ will be hit by residual waves. 
In three dimensions the situation is quite different. 
A flash of light at a point has an effect on the surface of a sphere around the point after a certain time interval but then no more.
There are no residual waves as those present on the water surface. 
The same is the case with sound waves, one has a pure propagation without residual waves; 
thus music can exist in $I\!\!R^{3}$.'' \cite{helg84} \ 
Waves fulfill a strong Huygens' Principle (HP) in our daily $I\!\!R^{3}$, however in $I\!\!R^{2}$ they do not. 
May be that such a ``global'' classification of waves in  $I\!\!R^{3}$, against  $I\!\!R^{2}$, is something misleading.
I this paper we meet the case that one solution of a given wave equation is a wake-free wave 
which individually fulfills a ``strong HP'', 
however, other solutions not.
It emerges the essential to define the use of the criterium ``HP'', either for a special solution, 
or for the solutions of all IVPs of a given wave equation.
Despite the widespread use of the term ``Huygens' principle'' in the physical community in the last centuries, 
its precise mathematical formulation was seldom clarified. 
In 1923, J.\,Hadamard formulated in mathematical terms  three different meanings 
of the HP he found in the literature of his time \cite{hadam}.
The general case is that the solution of the Cauchy IVP depends on the initial values 
in the full region inside the conoid of dependence;
however, a very special case emerges if the initial values are only used on the conoids surface area. 
This is the case with the strong validity of HP in the formulation 
of Hadamard's 'minor premise' \cite{hadam}, p.53\,f.
\newpage

{\bf Hadamard's syllogism} 
\small
\begin{itemize}
 \item \textbf{Major Premise}\\
        The action of phenomena produced
        at the instant $t=0$ on the state of matter at the later time
        $t=t_{0}$ takes place by the mediation of every intermediate
        instant $t=t'$, i.e. (assuming $0<t'<t_{0}$),  in order to
        find out what takes place for $t=t_{0}$, we can deduce from
        the state at $t=0$ the state at $t=t'$ and from the latter,
        the required state at $t=t_{0}$. 
 \item \textbf{Minor Premise}\\
        If, at the instant $t=0$ -- or more exactly
        throughout a short interval $-\epsilon\le t \le 0$
        -- we produce a luminous disturbance localized in the
        intermediate neighbourhood of $O$, the effect of it will be,
        for $t=t'$, localized in the immediate neighbourhood of the
        surface of the sphere with centre $O$ and radius $\omega t'$:
        that is, will be localized in a very thin spherical shell
        with centre $O$ including the aforesaid sphere.
 \item \textbf{Conclusion}\\
        In order to calculate the effect of our
        initial luminous phenomenon produced at $O$ at $t=0$, we may
        replace it by a proper system of disturbances taking place
        at $t=t'$ and distributed over the surface of the sphere with
        centre $O$ and radius $\omega t'$.
\end{itemize}
\normalsize

This paper is concerned with the second part, ``Hadamard's minor premise''.
It is known that the usual wave equation fulfills this strong HP in spaces of an odd dimension, $N=3,5,7,...$ only.
\\

\noindent
{\bf Definition} \ 
A wave equation where for every Cauchy data the solution fulfills the strong HP (the minor premise) is named a Huygens' equation.
\\

In our letter \cite{qubo13wII} we have shown that a special wave equation with a variable coefficient, with a special posed IVP, 
has a sharp solution being the contour hypersurfaces
of a given function in  $I\!\!R^{N+1}$, the level set wave equation (LSWE). 
The wake-free form of the contour hypersurfaces holds in every dimension $N$. 
The question emerges: are such LS wave equations also
Huygens' equations? In this paper, we obtain that the answer is no.
Different treatments are used for the result. In 
Section \ref{Kap3} we present a splitting formula for the wave operator which is important for the
dimension problem: our solution holds for odd as well as for even dimensions. 
The normal form of the LS wave equation is given in Section \ref{Kap4}, and 
in Section \ref{Kap5} we develop the curvature tensor and the Ricci tensor of the LSWE which we need in Section \ref{Kap6} to prove some
necessary conditions for a strong HP. 
In Sections \ref{Kap7} and \ref{Kap8} we calculate geodesic lines and a characteristic hypersurface of the LS wave equation,
and in Sections \ref{Kap9} and \ref{Kap10} we calculate the geodesic distance and a first part of the elementary solution, 
only to see that the criterium for such a first part for the validity of a strong HP is not fulfilled either. 
In a last Section we discuss the results:
it is not expectable that our LS wave equation fulfills the strong Huygens' principle.
%
\section{The level set wave equation (LSWE)} 
\label{Kap2}
%
We assume Cartesian coordinates ${\bf q}=(q^1,...,q^N)^T  \in I\!\!R^{N}$.
Let $\nu= V({\bf q})$ be the $N$-dimensional surface of interest in the space  $I\!\!R^{N+1}$ of 
points $({\bf q}^T,\nu)^T=(q^1,...,q^N,\nu)^T$. We assume that  $V({\bf q})$ is two times differentiable,
and ${\bf g}({\bf q})= \partial_{\bf q}V({\bf q})$ be its gradient vector, 
and ${\bf H}({\bf q})= \partial_{\bf q} {\bf g}^T = (\partial_{\bf q}\,\partial_{\bf q}^T\,V({\bf q}))$ 
be the matrix of its second derivatives, the Hessian. 
We form the scalar product of the gradient with itself. 
It produces the scalar function $\partial_{\bf q}V({\bf q})^T\,\partial_{\bf q}V({\bf q}) =: G({\bf q})$.
In the following, 
we treat regions without stationary points of the function $V({\bf q})$, thus there holds $G({\bf q}) > 0$.
Let us consider an $(N+1)$-dimensional linear wave operator with respect to the coordinates $q^1,...,q^N,\nu$, 
and with an additional term of the first order \cite{boqu12} 
\begin{equation}
 L\ :=  \Delta - G({\bf q})\ \frac{\partial^2}{\partial \nu^2} + Tr\,{\bf H}({\bf q})\ \frac{\partial}{\partial \nu}  \
\label{waveop}
\end{equation}
where $\Delta = \partial_{\bf q}^T \partial_{\bf q}=
(\partial/\partial q^1)^2+...+(\partial/ \partial q^N)^2$ is the 
Laplacian in $I\!\!R^N$, and $Tr\,{\bf H}$ is the sum of the diagonal entries of ${\bf H}$.
The operator $L$ is of the normal hyperbolic type by signature $++...+-$ while $G({\bf q})>0$, thus outside of 
stationary points of $V({\bf q})$. 
Note that $\nu$ is not the time variable.
We search for a solution of the equation   
\begin{equation} 
 L \psi({\bf q},\nu) = 0 \ .
\label{waveeq}
\end{equation}
The reasoning for the introduction of such an equation will come true with its success. 
\\

\noindent
{\bf Theorem} \\
Be $F$ a function of one real variable, $F: I\!\!R \rightarrow I\!\!R$, with first and second continuous derivations, then
\begin{equation}
\psi({\bf q},\nu) := F(V({\bf q}) - \nu)
\label{sol}
\end{equation} 
is a solution of the wave Eq.\,(\ref{waveeq}). 
\\

Such solutions are named progressing waves \cite{cour,fried}.
$V({\bf q})-\nu$ is the phase and $F$ is the wave form.
In the case of the classical wave equation with $N=3$ many special progressing waves have been known 
since the 18$^{th}$ century. 
Plane waves with phase $q^3 \pm c\,t$, $c=constant=1/G, \ t=\nu$, spherical waves 
with phase $R\pm c\,t$, $R=\sqrt{(q^1)^2+(q^2)^2+(q^3)^2}$ 
and an additional amplitude factor $1/R$.
Further solutions have been found in the 20$^{th}$  century named the Bateman-Hillion class. 
It is with phase $z+((q^1)^2+(q^2)^2)/w$ with $z=q^3 \pm c\,t$,
$w=q^3 \mp c\,t$ and amplitude factor $1/w$ \cite{hill93}. $w$ can still be extended to $w^*+C$ with a free constant $C$.
Generalizations of that class are further done by Borisov and Kiselev \cite{bor}.
These solutions have, comparably, simple and special phase functions: polynomials, logarithmic, trigonometric  or rational functions. 
The ansatz (\ref{sol}) is a very general case with a free function $V({\bf q})$. 
An early ansatz before was that of Gottlieb \cite{gott} with nonconstant coefficients, but for a circular symmetric equation.  
It was generalized by Bombelli et al.\,\cite{bom91}.
\\

\noindent
{\bf Proof of the Theorem} \\
The proof of (\ref{sol}) is straightforward: computing the differential expression $L \psi$, we get
\begin{equation}
L \psi = F''(..) \ [ \ (\partial_{\bf q} V)^T(\partial_{\bf q} V)  - G \ ] + 
         F'(..) \  [ \ \partial_{\bf q}^2 V - Tr\,{\bf H} \ ]  = 0 \ 
\label{proof}  
\end{equation}
per definition of $G$ and $Tr\,{\bf H}$.
$F$ can be an arbitrary function, because the coefficients of $F''$ and $F'$ are zero.
 $ \hfill{\rule{3mm}{3mm}}$ 
\\

 Of course, one can extend the usage of $F$ to distributions \cite{qubo13wII}.
\newpage

\begin{itemize}
\item[( i)]
 Now, we may assume that function G is prescribed and V is to determine. Then the vanishing of the  
first coefficient is the Hamilton-Jacobi equation or eikonal equation, 
a non-linear partial differential equation of the first order to search for V. 
It describes a relation between the level contours of a surface and its steepest descent lines \cite{bocre05}  
\begin{equation}
 (\partial_{\bf q} V({\bf q}))^T(\partial_{\bf q} V({\bf q})) = G({\bf q}) \ . 
\label{eikon}
\end{equation}
There is a large amount of methods to solve equations of this type \cite{osh88,set99,cec04}.
Here in Eq.(\ref{proof}) the eikonal is automatically fulfilled by the definition of $G({\bf q})$.
If one treats the limit $G({\bf q}) \rightarrow 0$, thus ${\bf q}$ would approach a stationary point of the surface $V({\bf q})$,
then also the eikonal (\ref{eikon}) degenerates.
Note that the eikonal also emerges in a variational theory of steepest descent lines \cite{bocre05}.
\item[(ii)]
We may also assume for the second coefficient that V is to determine. Then 
the second coefficient is named transport equation of the operator $L$ \cite{fried}, it emerges here in a simple version.
It is also automatically zero by the use of $Tr\,{\bf H}$.
\end{itemize}

\section{Dimension splitting}
\label{Kap3}
The solution (\ref{sol}) does not depend on the dimension of the coordinate space \cite{qubo13wII}. 
This is an astonishing result \cite{cour,fried,hill93,bor,gott,bom91}, 
and it is a hint that the strong Huygens' principle in the sense of Hadamard does not hold here \cite{hadam}.
We underline it with a partitioning of the coordinate space  $I\!\!R^{N}$ into two sets of variables
${\bf q}_{I}^T=(q_1,...,q_M)$ and  ${\bf q}_{II}^T=(q_{M+1},...,q_N)$ with $1 \le M<N$.
Due to the sum structure of the coefficients $G({\bf q})$ and $Tr\, {\bf H}({\bf q})$ we can split them into
\begin{equation}
G({\bf q})= G({\bf q}_I)+ G({\bf q}_{II}) {\rm \ \  and \ \ \ }  
\end{equation}
\begin{equation}
Tr\, {\bf H}({\bf q})= Tr\, {\bf H}({\bf q}_I)+ Tr\, {\bf H}({\bf q}_{II}) \ .
\end{equation}
Then we can write the differential eq.\,(\ref{waveeq}) as  
\begin{equation} 
 L \psi({\bf q},\nu) = \  L_I \psi({\bf q}_I,{\bf q}_{II},\nu) + L_{II} \psi({\bf q}_I,{\bf q}_{II},\nu) = 0 \ , 
\end{equation}
where the operator $L_I$ is
\begin{equation}
 L_I \ := \Big( \nabla_{{\bf q}_I}^2 - G({\bf q}_I)\ \frac{\partial^2}{\partial \nu^2} + 
               Tr\ {\bf H}({\bf q}_I)\ \frac{\partial}{\partial \nu} \Big) \
\label{waveI}
\end{equation}
acting in the  $I\!\!R^{M+1}$,  
and the same for $L_{II}$ acting in the  $I\!\!R^{N-M+1}$. 
If $\psi({\bf q}_I,{\bf q}_{II},\nu)$ is a solution of eq.\,(\ref{waveeq}) then also holds
$L_I \psi({\bf q}_I,{\bf q}_{II},\nu) =- L_{II} \psi({\bf q}_I,{\bf q}_{II},\nu) $.
The unique solution is  
$L_I \psi({\bf q}_I,{\bf q}_{II},\nu) = L_{II} \psi({\bf q}_I,{\bf q}_{II},\nu)=0$
which is independent of the partitioning. Especially, if $N$ was even, then a splitting into
$N=M+(N-M)$ with an odd $M$ produces two 'odd' operators.

We conclude: every ''full'' level hypersurface, $V({\bf q})=\nu$, of a surface, $V({\bf q})$, 
over coordinates in the $I\!\!R^{N}$ is an $(N-1)$-hypersurface and 
fulfills the $(N+1)$-dimensional differential eq.\,(\ref{waveeq}) with solution (\ref{sol}); 
but every restriction to subspaces of dimension $M<N$
forms a section in the subspace of the level hypersurface, of course on the same level $\nu$,
and it fulfills a reduced differential eq.\,(\ref{waveI}) in that subspace.
Thus, the operator $L$ in eq.\,(1) can be totally split into single summands
\begin{equation} 
 L \ = \sum_{i=1}^N \Big( \frac{\partial^2}{\partial q_i^2} -
 \Big( \frac{\partial V}{\partial q_i} \Big)^2\ \frac{\partial^2}{\partial \nu^2} + 
 \frac{\partial^2 V}{\partial q_i^2} \ \frac{\partial}{\partial \nu} \Big) \ =:  \sum_{i=1}^N L_i \ , 
\label{wavei} 
\end{equation}
and every part $L_i$ fulfills $L_i\psi =0$ where $\psi$ is a solution of $L \psi =0$. 
The proof (4) holds for every single $L_i$.
The sum (\ref{wavei}) is different from an operator sum in Ref.\,\cite{shi03}.
\section{Normal form of the LS wave equation}
\label{Kap4}
The part of the second order in operator (\ref{waveop}) applied to $\psi$ is:
\begin{equation}
\Delta \psi({\bf q},\nu) - G({\bf q}) \ \frac{\partial^2 }{ \partial^2 \nu} \psi({\bf q},\nu)
\label{waveeqP2} 
\end{equation} 
where $\psi ({\bf q},\nu)$ is
a field in a medium with ``slowness'' $G({\bf q})^{1/2}$ outside of stationary points of $V({\bf q})$. 
Note that the factor  $G({\bf q})$ depends only on the space variables {\bf q}, not on the variable $\nu$.
The contravariant metric matrix of the operator (\ref{waveeqP2}) is
\begin{equation} 
 (g^{I,J})_{I,J=1,...,N+1} = \left( \begin{array}{cc} {\bf E} & {\bf 0}  \\
     {\bf 0}^T & -G({\bf q}) \end{array}  \right) \  .
\end{equation}
${\bf E}$ is the N-dimensional unit matrix, and ${\bf 0}$ is the N-dimensional zero column vector. 
Because the metric matrix is only a diagonal matrix,
its inverse diagonal matrix is the covariant metric matrix
\begin{equation}
(g_{I,J})_{I,J=1,...,N+1} =
\left( \begin{array}{cc} {\bf E} & {\bf 0}  \\
     {\bf 0}^T & -1/G({\bf q}) \end{array}  \right) \  .
\label{metric}
\end{equation} 
%
%
Let the positive number $\gamma$ be the absolute value of its determinant, here $\gamma =1/G({\bf q})$.
\ \
The $\nu$-part of the operator (\ref{waveeqP2}) can be written
\begin{equation} 
-G({\bf q})^{1/2}  \frac{\partial }{ \partial \nu} \left( G({\bf q})^{-1/2} G({\bf q}) \frac{\partial }{ \partial \nu}  \right)\psi({\bf q},\nu) =
\frac{1}{\sqrt{\gamma}} \frac{\partial }{ \partial \nu} \left( \sqrt{\gamma} g^{N+1,N+1} \frac{\partial }{ \partial \nu}  \right)\psi({\bf q},\nu) \ .
\label{old9}
\end{equation} 
For the space variables we have
\begin{equation}
G({\bf q})= \sum_{i=1}^{N} \left( \frac{\partial }{ \partial q^i} V({\bf q}) \right)^2, {\rm \  thus \ \ } 
 \frac{\partial }{ \partial q^k} G({\bf q})= 2 \sum_{i=1}^{N} V_{q^i}({\bf q}) V_{q^iq^k}({\bf q}) =: 2\ {\cal H}_k({\bf q})\  .
\label{calH}
\end{equation} 
Here, ${\cal H}_k$ is the Euclidean scalar product of the gradient of $V$ with the $k$-th column of the Hessian matrix.
The coefficients are $g^{i,j}= \delta ^{i,j}$ \ for $i,j=1,...,N$, being the metric fundamental tensor of the  $I\!\!R^N$.
For the summands of the wave operator of Eq.(\ref{waveeqP2}) we treat the ansatz 
\begin{equation}
\frac{1}{\sqrt{\gamma }}\ \frac{\partial }{\partial q^i} \left( \sqrt{\gamma} \frac{\partial }{ \partial q^i}  \right)\psi({\bf q},\nu) \ . 
\end{equation} 
If we differentiate the factor $\sqrt{\gamma}$ inside the formula, we get
\begin{equation}
\frac{\partial ^2}{ \partial q^{i\,2}} \psi({\bf q},\nu) + G({\bf q})^{1/2} \ (-\frac{1}{2} )\frac{1}{ G({\bf q})^{3/2}} \ 2  {\cal H}_i\ 
\frac{\partial }{\partial q^i }\psi({\bf q},\nu)  \ 
\end{equation} 
with the product rule and with (\ref{calH}).
Consequently, we get for the operator (\ref{waveeqP2}) 
\begin{equation}
\Delta \psi({\bf q},\nu) - G({\bf q}) \ \frac{\partial^2 }{ \partial \nu\,^2} \psi({\bf q},\nu) = 
\label{nf}
\end{equation}
$$
 \frac{1}{\sqrt{\gamma }}\ \frac{\partial }{\partial q^I} \left( \sqrt{\gamma}  g^{I,J} 
 \frac{\partial }{ \partial q^J}  \right)\psi({\bf q},\nu)
+ \frac{1}{G({\bf q})} \sum_i {\cal H}_i({\bf q})\ \frac{\partial }{\partial q^i}\psi({\bf q},\nu) \ , 
\nonumber
$$
and the first part of the right hand side of Eq.(\ref{nf}) 
builds the normal form of a second order wave equation with (N+1) variables and with variable coefficients.
We will use the convention that an index $j$ is to be summed over 1 to $N$ whenever it is repeated in a single term, as well as a $J$ is to be summed over 1 to $N+1$. 
This has to to be done in the parts of the right hand side of Eq.(\ref{nf}).
Note that the first order part of operator (\ref{waveop}) becomes
\begin{equation}
\frac{1}{G} \sum_i^N {\cal H}_i \frac{\partial }{\partial q^i}  + Tr\,{\bf H} \frac{\partial }{\partial \nu } \ ,
\label{cov}
\end{equation} 
but a zeroth order part is missing. The expressions ${\cal H}_i$ are defined in Eq.(\ref{calH}).
%
\section{Curvature tensor and Ricci tensor of the metric connected with the LS wave equation}
\label{Kap5}
%
%
For the general formulae of the following see any text book of differential geometry. 
With  the metric matrix (\ref{metric}) the corresponding coordinate space becomes an (N+1)-dimensional manifold. 
Its metric is given by $g_{II}\,dq^I\,dq^I$ in the used coordinates.
We can calculate the components of the affine connection (Christoffel symbols of second kind) determined by
\begin{equation}
\Gamma_{JK}^I =\frac{1}{2} g^{IL} \left( \frac{\partial g_{LJ}}{ \partial q^K} +
\frac{\partial g_{KL}}{ \partial q^J} -  \frac{\partial g_{JK}}{ \partial q^L}  \right)
\label{eq1}
\end{equation}
where 
$1 \le I,J,K,L \le N+1 $. In the following, we will use lower-case letters $1 \le i,j,k \le N $, 
and we often separately treat the (N+1)th coordinate, $\nu$.
Because the metric is somewhat simple, we get most of the $\Gamma$ components to be zero.
Nonzero are only the three kinds of the following components:
\begin{equation}
\Gamma_{N+1\,N+1}^{\, i} =\frac{1}{2}  \left(0+0 -  \frac{\partial g_{N+1\,N+1}}{ \partial q^i}  \right) =
\frac{1}{2}  \left( - \frac{\partial }{ \partial q^i} \frac{-1}{G} \right) =
\frac{-1}{G^2} \,\sum_k^N   \frac{\partial V}{ \partial q^k} \, \frac{\partial^2 V}{ \partial q^k  \partial q^i} \ ,
\label{eq2}
\end{equation}
where we use $ G= \sum_k^N (\frac{\partial V}{ \partial q^k})^2 $. Thus, we have by definition (\ref{calH})
\begin{equation}
\Gamma_{N+1\,N+1}^{\, i} =  \frac{-1}{G^2}  {\cal H}_i \ .
\label{eq3}
\end{equation}
 Analogously, we find
\begin{eqnarray}
\Gamma_{N+1\,k}^{\, N+1} =\frac{1}{2}( -G ) \left( \frac{\partial g_{N+1\,N+1}}{ \partial q^k} +0 -0  \right) =
\frac{1}{2} ( -G ) \left(\frac{\partial }{ \partial q^k} \frac{-1}{G} \right)
\nonumber
\\ 
= \frac{-1}{G } \,\sum_i^N   \frac{\partial V}{ \partial q^i} \, \frac{\partial^2 V}{ \partial q^i  \partial q^k} \ . \quad \quad \quad \quad \quad \quad \quad \quad 
\label{eq4}
\end{eqnarray}
Thus, 
\begin{equation}
\Gamma_{N+1\,k}^{\, N+1} = \frac{-1}{G} {\cal H}_k \  {\rm \ \quad \ \quad \ \quad  and \ analogously \ \ \quad \ \quad  \ \quad }
\Gamma_{j\ N+1}^{\, N+1} = \frac{-1}{G} {\cal H}_j \ .
\label{eq5}
\end{equation}

\noindent
Using the Christoffels we can also develop the operator (\ref{nf}) by the formula
\begin{equation} 
 g^{I\,J}\, \nabla_I \nabla_J =
\nabla^J \nabla_J =  g^{I\,J}\, \partial_I \partial_J -  g^{I\,J}\,\Gamma_{I\,J}^K \partial_K \ .
\end{equation}

\noindent
We define new helpful symbols \cite{schi} using the Christoffels  (\ref{eq3}) and (\ref{eq5})
\begin{equation} 
\Gamma_{I\,J\,K}^{\quad\quad L} := \frac{\partial}{\partial q^I}  \Gamma_{J\,K}^{\,L} - \sum_M  \Gamma_{I\,K}^{\,M}\, \Gamma_{J\,M}^{\,L} \ .
\label{symb}
\end{equation}
Because the metric matrix (\ref{metric}) is very simple, we get only some single $\Gamma_{I\,J\,K}^{\quad\quad L}$ not zero:
\begin{eqnarray}
\Gamma_{i\,j\,N+1}^{\quad\quad\quad N+1}\ =\quad \frac{\partial}{\partial q^i} \Gamma_{j\,N+1}^{\,N+1} - \Gamma_{i\,N+1}^{\,N+1}\, \Gamma_{j\,N+1}^{\,N+1} \
=\quad \frac{\partial}{\partial q^i}( \frac{-{\cal H}_j}{G} ) - \frac{1}{G^2} {\cal H}_i {\cal H}_j
\nonumber
\\
=\quad \frac{1}{2} \frac{\partial}{\partial q^i} \left( G \frac{\partial}{\partial q^j} \frac{1}{G}\right) -
\frac{1}{4} G^2 (\frac{\partial}{\partial q^i} \frac{1}{G}) 
(\frac{\partial}{\partial q^j} \frac{1}{G})  \quad \quad \quad \quad
\label{symb1}
\end{eqnarray}
where only the summand with $M=N+1$ gives a contribution, and
\begin{eqnarray}
\Gamma_{i\,N+1\,N+1}^{\quad\quad\quad\quad l} = 
\quad \frac{\partial}{\partial q^i} \Gamma_{N+1\,N+1}^{\,l} - \Gamma_{i\,N+1}^{\,N+1}\, \Gamma_{N+1\,N+1}^{\,l} \
=\quad \frac{\partial}{\partial q^i}( \frac{-{\cal H}_l}{G^2} ) - \frac{1}{G^3} {\cal H}_i {\cal H}_l 
\nonumber
\\
=\quad \frac{1}{2} \frac{\partial^2}{\partial q^i \partial q^l} \frac{1}{G} 
- \frac{1}{4}\, G\,  (\frac{\partial}{\partial q^i} \frac{1}{G}) 
(\frac{\partial}{\partial q^l} \frac{1}{G}) \ , \quad\quad\quad\quad
\label{symb2} 
\end{eqnarray}
but it is 
\begin{equation}
\Gamma_{N+1\,j\,N+1}^{\quad\quad\quad\quad l} = 0 \ , 
\label{symb20} 
\end{equation}
which we need below. It is further
\begin{equation}
\Gamma_{N+1\,N+1\,k}^{\quad\quad\quad\quad l} = 
\frac{-G}{4} (\frac{\partial}{\partial q^k}\frac{1}{G})
             (\frac{\partial}{\partial q^l}\frac{1}{G}) \ ,
\label{symb26n}
\end{equation}
\begin{equation}
\Gamma_{N+1\,j\,k}^{\quad\quad \ N+1} = 
\frac{-G^2}{4} (\frac{\partial}{\partial q^j}\frac{1}{G})
               (\frac{\partial}{\partial q^k}\frac{1}{G}) \ ,
\label{symb22a}
\end{equation}
\begin{equation}
\Gamma_{i\,N+1\,k}^{\quad\quad N+1} = \quad \frac{\partial}{\partial q^i}(\frac{-{\cal H}_k}{G} ) 
= \frac{1}{2} \frac{\partial}{\partial q^i}(G \frac{\partial}{\partial q^k} \frac{1}{G} ) \ ,
\label{symb22}
\end{equation}
and 
\begin{equation}
\Gamma_{N+1\,N+1\,N+1}^{\quad\quad\quad\quad \ N+1} = - \Gamma_{N+1\,N+1}^{\,m}\, \Gamma_{N+1\,m}^{\,N+1} 
=- \frac{1}{4}\, G\, \sum_m^N  (\frac{\partial}{\partial q^m} \frac{1}{G})^2
= \frac{-1}{G^3} \sum_m {\cal H}_m^2 \ .
\label{symb3}
\end{equation}

\noindent
By an antisymmetrisation in I,J of symbols (\ref{symb}) the Riemannian curvature tensor  emerges \cite{gall}, p.106
\begin{equation}
  R_{I\,J\,K}^{\quad\quad\,L} =  \Gamma_{I\,J\,K}^{\quad\quad\,L} -
  \Gamma_{J\,I\,K}^{\quad\quad\,L} {\rm \quad with \quad} R_{I\,J\,K}^{\quad\quad\,L} =-R_{J\,I\,K}^{\quad\quad\,L} \ .
\end{equation}
With the symbol (\ref{symb1}) being symmetric in $I,J$ we get $R_{i\,j\,N+1}^{\quad\ \ N+1}=0$.
With Eq.(\ref{symb26n}) we get $R_{N+1\,N+1\,k}^{\quad\quad\quad\quad l} = 0$, 
and analogously with (\ref{symb3}) we have $R_{N+1\,N+1\,N+1}^{\quad\quad\quad\quad \ N+1}=0$.
But it is with (\ref{symb2}) and (\ref{symb20})
\begin{equation}
R_{i\,N+1\,N+1}^{\quad\quad\quad \  l} = 
\Gamma_{i\,N+1\,N+1}^{\quad\quad\quad \ l} \ \ne 0 \  , 
\label{cur2}
\end{equation}
and with (\ref{symb22a}) and (\ref{symb22})
$$
R_{N+1\,j\,k}^{\quad\quad N+1}= 
\Gamma_{N+1\,j\,k}^{\quad\quad N+1} \ - \Gamma_{j\,N+1\,k}^{\quad\quad N+1}  = -R_{j\,N+1\,k}^{\quad\quad N+1}
$$
\begin{equation}
= -\frac{1}{2} \frac{\partial}{\partial q^j}(G \frac{\partial}{\partial q^k} \frac{1}{G} ) \ 
-  \frac{G^2}{4} (\frac{\partial}{\partial q^j}\frac{1}{G})
   ( \frac{\partial}{\partial q^k} \frac{1}{G})  \ .
\label{cur22}
\end{equation}

\noindent
The Ricci tensor is defined by
\begin{equation}
R_{J\,K} = \sum_M R_{M\,J\,K}^{\quad\quad\,M}  {\rm \quad with \quad}R_{J\,K} = R_{K\,J} \ ,
\end{equation}
and the nonzero Ricci components will be obtained by a summation over $i=l=m$ of the components of the kind (\ref{cur2}) giving
\begin{equation}
R_{N+1\,N+1} = \sum_m^N \left( 
 \frac{1}{2} \frac{\partial^2 }{\partial q^{m\,2}} \frac{1}{G}
- \frac{1}{4}\, G\, (\frac{\partial}{\partial q^m} \frac{1}{G})^2 \right)  \ ,
\label{ricci}
\end{equation}
or by (\ref{cur22}) with only one summand, giving
\begin{equation}
 R_{j\,k} = 
-\frac{1}{2} \frac{\partial }{\partial q^{j}} (G \frac{\partial}{\partial q^k}  \frac{1}{G})
- \frac{1}{4}\, G^2\, (\frac{\partial}{\partial q^j} \frac{1}{G})
\, (\frac{\partial}{\partial q^k} \frac{1}{G})
\ .
\label{ricci2}
\end{equation} 
They are not zero, in general, while the surface V is not flat.

\noindent
The curvature scalar is
$$
R = g^{J\,K}\, R_{J\,K} \  = \sum_{j\,k}^N \delta^{j\,k} R_{j\,k} - G R_{N+1\,N+1} 
$$
\begin{equation}
= \sum_j \left( \frac{-3}{2G^2} \ ( \frac{\partial}{\partial q^{j}} G )^2 +
  \frac{1}{G} \ ( \frac{\partial^2}{\partial q^{j\,2}} G )  \right) 
= \sum_j \left( \frac{-6\,{\cal H}_j^2}{G^2} \  +
  \frac{2}{G} \ ( \frac{\partial}{\partial q^j} {\cal H}_j )  \right) \ . 
\label{curvR}
\end{equation}  
%
\section{Negative results concerning conditions for a strong Huygens' Principle}
\label{Kap6}
The general wave equation is written in a coordinate invariant form
\begin{equation}
 L\  \Psi \ =  g^{I\,J} \, \nabla_I \, \nabla_J \Psi + A^I \, \nabla_I \Psi + C \, \Psi = 0 \ .
\label{invari}
\end{equation}
In Eq.(\ref{waveop}) we have the case $C=0$, and coefficient $A^{N+1}=Tr {\bf H}$.
With the covariant form of the wave equation (\ref{cov}), we have additionally the contravariant 
coefficients  $A^{i}= \frac{1}{G}{\cal H}_i $.

The proof for the validity of the strong HP property is quite difficult since the validity conditions involve the coefficients of 
operator (\ref{nf}) in a very indirect and complicated manner.
We calculate a first condition for HP given in Ref.\,\cite{mc}, Eq.(1.19), or see also Ref.\,\cite{sch78}.
It is with n=N+1 here (if n=4, then -1/6 is the correct factor for R)
\begin{equation}
C -\frac{1}{2} \sum_I^{N+1} \nabla_I A^I -\frac{1}{4}\,\sum_I^{N+1} A_I\, A^I - \frac{N-1}{4\,N}\, R = 0 \ .
\label{mcL}
\end{equation}

The covariant coefficients for (\ref{mcL}) are $A_{i}= \delta_{i\,j} A^{j} = A^{i} = \frac{1}{G}{\cal H}_i$ for $i=1,...,N$, 
and $A_{N+1}= g_{N+1\,N+1} A^{N+1} =-\frac{1}{G} Tr {\bf H}$.
The symbol $\nabla_I$ means covariant differentiation
\begin{equation}
 \nabla_K A^I  = \partial_K A^I +  \sum_J^{N+1} \Gamma_{K\,J}^I A^J \ ,
\end{equation}
in words: 
the covariant differentiation is the usual derivative along the coordinates with correction terms which tell us how the coordinates change themselves.
It is here
\begin{equation}
 \nabla_k A^i  = \partial_k A^i \ ,
\end{equation}
because the corrections are $\Gamma_{k\,J}^i=0$ throughout. Thus, the $q^k$ do not depend on all other coordinates.
However, it is 
\begin{equation}
 \nabla_{N+1} A^{N+1}  = 0 + \sum_J^{N+1} \Gamma_{N+1\,J}^{N+1} A^J = \sum_j^{N} \Gamma_{N+1\,j}^{N+1} A^j \ .
\end{equation}
Though the coefficient $A^{N+1}$ does not depend on $\nu$, its covariant derivation to $\nu$ is corrected by terms which include the other $A^{j}$  here.
A change along $\nu$ depends on all other coordinates $q^i$ in a highly nonlinear kind. 
This part in Eq.(\ref{mcL}) becomes
\begin{eqnarray}
-\frac{1}{2} \sum_I^{N+1} \nabla_I A^I = -\frac{1}{2} \sum_i^N \frac{\partial}{\partial q^i}(\frac{1}{G}{\cal H}_i) + 0 
-\frac{1}{2} \sum_j^N \Gamma_{N+1\,j}^{N+1}\, A^j
=
\nonumber
\\
=   \frac{1}{G^2}\, \sum_i^N {\cal H}_i^2- \frac{1}{2 G}\,\sum_i^N 
(\frac{\partial}{\partial q^i}{\cal H}_i) + \frac{1}{2 G^2}\, \sum_j^N {\cal H}_j^2 \ .
\label{covI}
\end{eqnarray}
The next summand in (\ref{mcL}) is 
\begin{equation}
  -\frac{1}{4} \sum_I^{N+1} A_I \, A^I =
 -\frac{1}{4 G^2}\sum_i^N  {\cal H}_i^2 +\frac{1}{4 G} (Tr {\bf H})^2 \ .
\label{sumI}
\end{equation}
All in all,
the condition (\ref {mcL}) has the following summands, if we use (\ref{curvR}), (\ref{covI}) and (\ref{sumI}) and sort equal terms
\begin{eqnarray}
 \frac{1}{4\,G}\, (Tr {\bf H})^2   - \frac{1}{2\,G} \sum_i^N 
 \frac{\partial}{\partial q^{i}} {\cal H}_i \ \left(1 +\frac{N-1}{N} \right)
+\frac{1}{4\,G^2}  \sum_i^N  {\cal H}_i^2  \ \left( 5+6\,\frac{N-1}{N} \right) \ .  
\label {mcL2}
\end{eqnarray}
The expression seems never to be zero.
\newpage 

%
\noindent
{\bf Conclusion } \ 
The first condition is not fulfilled for the validity of a strong HP of the wave operator (\ref{waveop}), 
for Hadamard's 'minor premise'. Thus, wave operator (\ref{waveop}) is not a Huygens' operator.
\\ 

A second  condition  only concerns the coefficients $A^I$ and the metric \cite{mc}.
We define (by another capital letter than Ref.\,\cite{mc} to differentiate from the Hessian)
\begin{equation}
 K_{I\,J} := \frac{1}{2} \left( \frac{\partial}{\partial q^J} A_I -  \frac{\partial}{\partial q^I} A_J \right) \ .
\label{c2}
\end{equation}
Because of symmetries, it holds $K_{i\,j} =0$ for $i,j=1,...,N$, as well as $K_{N+1\,N+1} =0$. 
Only two kinds of components are not zero:
\begin{equation}
  K_{N+1\,j} = \frac{1}{2} \frac{\partial}{\partial q^j} A_{N+1} = - K_{j\,N+1} = 
  -\frac{1}{2} \frac{\partial}{\partial q^j}   \frac{Tr {\bf H}}{G} \ .
\end{equation}
The second condition of \cite{mc} is
\begin{equation}
\nabla^J K_{I\,J} = 0 \ .
\end{equation}
It is defined \cite{fried}
\begin{equation}
 \nabla_L A_I  = \partial_L A_I - \sum_M^{N+1} \Gamma_{L\,I}^M A_M \ ,
\end{equation}
and \cite{will}
\begin{equation}
 \nabla_L K_{I\,J}  = \partial_L K_{I\,J} - \sum_M^{N+1} \Gamma_{L\,I}^M K_{M\,J} 
                                          - \sum_M^{N+1} \Gamma_{L\,J}^M K_{I\,M} \ ,
\end{equation}
thus we get
\begin{equation}
\nabla^J K_{I\,J} = g^{JL} \nabla_L K_{I\,J}=  g^{JL} 
\left( \partial_L K_{I\,J} - \sum_M^{N+1} \Gamma_{L\,I}^M K_{M\,J} - \sum_M^{N+1} \Gamma_{L\,J}^M K_{I\,M} 
\right) \ . 
\end{equation}
If
$1\le i \le N$, then 
\begin{equation}
\nabla^J K_{i\,J} = g^{JL} \nabla_L K_{i\,J}= g^{N+1\,L} \partial_L K_{i\,N+1} 
- g^{j\,N+1} \Gamma_{N+1\,i}^{N+1} K_{N+1\,j} -  g^{JL} \Gamma_{L\,J}^{N+1} K_{i\,N+1}=0  \ . 
\end{equation}
However, if $I=N+1$ it is
$$
\nabla^J K_{N+1\,J} = g^{JL} \nabla_L K_{N+1\,J}=  g^{jL} \partial_L K_{N+1\,j} -g^{j\,l} \Gamma_{l\,N+1}^{N+1} K_{N+1\,j} 
$$
$$
- g^{N+1\,N+1} \sum_m^{N} \Gamma_{N+1\,N+1}^m ( K_{m\,N+1} + K_{N+1\,m} ) 
$$
\begin{equation}
= \sum_j^N \partial_j K_{N+1\,j}  -\sum_j^N \Gamma_{j\,N+1}^{N+1} K_{N+1\,j} 
= \sum_j^N  \frac{1}{2}\left( \frac{\partial^2}{\partial q^{j\,2}} 
+ \frac{1}{G} {\cal H}_j \frac{\partial}{\partial q^{j}} \right) \frac{-Tr\,{\bf H}}{G} \ne 0 \ ,
\end{equation}
in the general case. So also the second condition of Ref.\,\cite{mc} for a strong HP is not fulfilled; 
and we stop these calculations here.
%
\section{Special geodesics for the LS wave equation}
\label{Kap7}
Generally, the geodesics are determined by a system of ordinary differential equations:  
\begin{equation}
\frac{d^2 q^I}{ds^2} + \Gamma_{J\,K}^{\, I} \frac{d q^J}{ds} \frac{d q^K}{ds} =0, \ \ I,J,K=1,...,N+1 \ ,
\label{eq6}
\end{equation}
where $s$ is the arclength parameter (see a remark below, after Eq.(\ref{eq17})), and it is summarized over J and K. 
Again, we use $q^{N+1}=\nu $, and the lower-case indexes for the $q^i$, and get the (N+1) coupled equations:
\begin{eqnarray}
\frac{d^2 \nu}{ds^2} - \frac{2}{G}\, \frac{d \nu}{ds} ( \sum_j^N \frac{d q^j}{ds}  {\cal H}_j ) & =0, &
\label{eq7}
\\
\frac{d^2 q^i}{ds^2} - \frac{{\cal H}_i}{G^2}\,  (\frac{d \nu}{ds})^2   & =0, & \ \  i=1,...,N \ .
\label{eq8}
\end{eqnarray}
Eq.\,(\ref{eq7}) can be simplified 
\begin{equation}
\frac{d^2 \nu}{ds^2} - \left( \frac{d}{ds} log\,(G) \right) \frac{d \nu}{ds} =0 \ .
\label{eq9} 
\end{equation}
The equation has a solution 
\begin{equation}
 \nu (s) = \int_{s_o}^s G({\bf q}(\tau) ) d\tau  {\rm \ \ \ with \ \ \ }  \frac{d \nu}{ds}(s) =  G({\bf q}(s)) \ .
\label{eq10}
\end{equation}
The function  $G({\bf q})$ determines the relation between the level, $\nu$, and the arclength parameter, $s$.
Using it in Eq.(\ref{eq8}), we can then try the steepest ascent ansatz for the (formally decoupled) solution for the space coordinates by
\begin{equation}
 \frac{d q^i}{ds}(s)=  \frac{\partial V}{ \partial q^i}({\bf q}(s)) \ , \ \  i=1,...,N  \ .
\label{eq11}
\end{equation}
For the proof of Eq.(\ref{eq11}) we do a next derivation to $s$. It is
\begin{equation}
\frac{d^2 q^i}{ds^2} =
 \,\sum_k^N \frac{\partial }{ \partial q^k} \, (\frac{\partial V}{ \partial q^i}) \  \frac{d q^k}{ds} =  
 \,\sum_k^N \frac{\partial V}{\partial q^k} \, \frac{\partial^2 V}{ \partial q^k  \partial q^i} = {\cal H}_i \ ,
\label{eq11-2}
\end{equation}
thus, the correct expression in Eq.(\ref{eq8}). 
Of course, all coordinates depend on each other because they all are used in $G({\bf q})$ and $V({\bf q})$. 
Thus indeed, we get the steepest ascent lines for the surface, V({\bf q}), for being special geodesics. 
And $\nu$ develops corresponding to Eq.(\ref{eq10}).

Note: another, trivial solution of Eq.\,(\ref{eq7}), in contrast to Eq.(\ref{eq10}), is $d\nu/ds=0$ which 
allows straight lines in the $q^i$-space for the geodesics, in a hyperplane $\nu=$constant.
\section{The characteristic surface and the initial value problem (IVP)}
\label{Kap8}
The method of characteristics is a technique developed for first order equations, though it is 
also valid for a second order hyperbolic equation. 
The idea is to reduce a partial differential equation to a family of ordinary differential equations 
along which the solution can be integrated from some initial data given on a suitable hypersurface.
First, we look for the characteristic manifold of operator (\ref{waveop}) \cite{cour}, part II, Chap.VI, \S 2.
We can treat the direct solution of Eq.\,(\ref{waveeq})
\ $ \psi({\bf q},\nu) = F(V({\bf q})-\nu) $.  
Then we have the phase function
\begin{equation}
S({\bf q},\nu) = V({\bf q}) - \nu
\label{phase}
\end{equation} 
which directly fulfills the characteristic equation \cite{cour} pertaining to operator $L$ of Eq.\,(\ref{waveop})
\begin{equation}
 (\nabla_{\bf q} S({\bf q},\nu ))^T(\nabla_{\bf q} S({\bf q},\nu )) - 
 G({\bf q})\left( \frac{\partial}{\partial \nu} S({\bf q},\nu )\right)^2  = 0\ . 
\label{char}
\end{equation}
It is identical with the eikonal Eq.(\ref{eikon}). 
On the other hand, if we treat a steepest ascent trajectory $({\bf q}(s),\nu(s))$ then $S({\bf q},\nu)$ is constant along this line. It is 
\begin{equation}
 \frac{d}{ds} S({\bf q}(s),\nu(s)) =\frac{\partial V}{\partial q^i}  \frac{d\,q^i}{ds}- \frac{d\,\nu}{ds}
=  \sum_{i=1}^N  ( \frac{\partial}{\partial q^i} V({\bf q}) )^2 - G({\bf q})=0 , 
\end{equation}
where we used the property of steepest ascent lines. The constant is zero, of course, 
because along the steepest ascent line it holds
$V({\bf q}) = \nu$. 
\\

On the other hand, with Friedlander \cite{fried}, p.79, a characteristic surface, $\varphi ({\bf q},\nu)=0 $,
is a null-surface in the metric of the space used, thus (\ref{metric}). The coefficient of the $\nu$ coordinate is not constant, 
thus we cannot treat distances over global regions by the usual 'Pythagoras'.
We have to use local neighbourhoods and differentials.
With the local point of view, to formally find the characteristic manifold, 
we have to treat the quadratic form which is connected with the metric (\ref{metric})
\begin{equation}
 \sum_{i=1}^N  ( \frac{\partial}{\partial q^i} \varphi ({\bf q},\nu))^2 - 
\frac{1}{G({\bf q})} (\frac{\partial}{\partial \nu} \varphi ({\bf q},\nu))^2=0 \ .
\label{eqx4}
\end{equation}
It describes a cone for every point $({\bf q},\nu)$.
If we assume a curve embedded in the characteristic surface, $({\bf q}(s),\nu(s))$, $s$ is the arclength, 
then the tangent direction has to be a null vector
\begin{equation}
 \sum_{i=1}^N  ( \frac{d q^i}{d s})^2 - 
\frac{1}{G({\bf q})} (\frac{d \nu}{d s})^2=0 \ .
\label{tgx4}
\end{equation}
Using the steepest ascent equation with eqs.(\ref{eq11}) and (\ref{eq10}) which is a null geodesic,
we again get Eq.\,(\ref{char}) by 
\begin{equation}
 \sum_{i=1}^N  ( \frac{\partial}{\partial q^i} V({\bf q}) )^2 - 
\frac{1}{G({\bf q})} (G({\bf q}))^2=0 \ .
\label{eqx44}
\end{equation}
Note, by the use of the steepest ascent equation in the ansatz of the characteristics, we only use one sort of ascending lines,
and thus exclude the full 'conoid'. 
The result is then, indeed, the function V({\bf q}) only. It emerges as an envelope of the current cones.

We find that $S({\bf q},\nu) = \varphi ({\bf q},\nu)$ is a special characteristic manifold. 
The ''space-like'' derivatives of function S are the derivatives of the surface, V({\bf q}). 
Consequently, V({\bf q}) itself is a characteristic surface to operator (\ref{waveop}).
\ 
Characteristic surfaces $S({\bf q},\nu) =0$ are considered as potential carriers of wave fronts \cite{cour}.
But the hypersurfaces $S({\bf q},\nu) =0$ 
with fixed $\nu$ are also the contours of the surface, V.

A further trivial solution of system (\ref{char}) is the function $S_+({\bf q},\nu) = V({\bf q}) + (\nu -2\nu_0 )$ with any fixed value $\nu_0$.
Generally, it is not of interest for us; we are interested in the contours $S=0$ of the phase function (\ref{phase}) only. 

However,
using together phase functions like $S$ and $S_+$ we can pose an initial value problem for Eq.(\ref{waveeq}) \cite{qubo13wII}.
We assume two functions $F$ and $D$ of one real variable, $F,D: I\!\!R \rightarrow I\!\!R$, with some continuous derivations, and we put 
\begin{equation}
\psi({\bf q},\nu) =\frac{1}{2} \left\{ F(V({\bf q}) + \nu)+ F(V({\bf q}) - \nu) +  
                                       D(V({\bf q}) + \nu)- D(V({\bf q}) - \nu) \right\} \ .
\label{ivp}
\end{equation} 
It is a solution like (\ref{sol}) of the wave Eq.\,(\ref{waveeq}), as well. 
We get back (\ref{sol}) if we use $D=F$. 
Solution (\ref{ivp}) fulfills the following initial values, for example for $\nu_o = 0$
\begin{equation}
\psi({\bf q},\nu_o) = F(V({\bf q}))
\label{ivp_0}
\end{equation}
and
\begin{equation}
\frac{\partial}{\partial \nu} \psi({\bf q},\nu)|_{\nu = \nu_0} = D'(V({\bf q})) \ .
\label{ivp_1}
\end{equation}
Of course, the free functions $F(x)$ and $D(x)$ are here only one-dimensional, 
because the solution (\ref{ivp}) is adapted to the $N$-dimensional surface, $V({\bf q})$, 
being the basic of the operator (\ref{waveeq}).
Thus, (\ref{ivp_0}) and (\ref{ivp_1}) do not pose the most general initial value problem.
However, they are for (\ref{waveeq}) a well posed Cauchy initial value problem.  
%
\section{The geodesic distance}
\label{Kap9}
Note that the considerations in this Section are again of local character.
\\

\noindent
{\bf Lemma}\\
Be $({\bf q}(s)^T,\nu(s))$ a point on a general geodesic line. It holds
\begin{equation}
\frac{d}{ds} \left[ \sum_{k=1}^N (\frac{d q^k}{ds})^2 -\frac{1}{G} (\frac{d \nu}{ds})^2 \right] =0 \ .
\label{eq12}
\end{equation}
It is an important relation which enfolds somewhat of the character of a geodesic distance.

\noindent
{\bf Proof} \\ 
If one derivates, it holds
\begin{equation}
= \sum_{k=1}^N 2 (\frac{d q^k}{ds}) (\frac{d^2 q^k}{ds^2}) -(\frac{d}{ds} \frac{1}{G}) (\frac{d \nu}{ds})^2 
- \frac{2}{G} (\frac{d \nu}{ds})  (\frac{d^2 \nu}{ds^2}) 
\nonumber
\end{equation}
and with the differential equations of the general geodesics (\ref{eq7}),(\ref{eq8}) it is
\begin{equation}
= \left[ -\frac{d}{ds} \frac{1}{G}  -\frac{d}{ds} \frac{1}{G}  -\frac{2}{G^2} \, \frac{d}{ds} G \right]\, (\frac{d^2 \nu}{ds^2}) \ .
\nonumber
\end{equation}
But the last summand is $+2  \frac{d }{ds}  \frac{1}{G} $. Thus, the given relation is correct, cf.\,also Ref.\,\cite{bocre05}.
A consequence of Eq.\,(\ref{eq12}) is
\begin{equation} 
\sum_{k=1}^N (\frac{d q^k}{ds})^2 -\frac{1}{G} (\frac{d \nu}{ds})^2  = const = D \ ,
\label{eq15}
\end{equation}
where $D$ may differ with different geodesic lines, but it is constant along a given geodesic.
If we integrate the square root of the equation along a geodesic line ${\bf q}(\tau)$ it is with $J,K=1,...,N+1$ in the general case  
\begin{equation}
 \int_o^s \left[ g_{JK}({\bf q}(\tau)) \, \dot{q}^J(\tau)\, \dot{q}^K(\tau)  \right]^{1/2} d\tau  \  ,
\label{eq16}
\end{equation}
thus here
\begin{equation}
{\cal J}({\bf q},\nu ,s)= \int_o^s \left[ \sum_{k}^N \dot{q}^k(\tau)^2  -\frac{1}{G} \dot{\nu}(\tau)^2 \right]^{1/2} d\tau \ = \ D^{1/2 } \, s \ ,
\label{eq17}
\end{equation}
Eq.(\ref{eq17}) is the formula for the arclength of the geodesics between a point $({\bf q}_o^T,\nu_o)$ at parameter $s=0$, 
and the final point $({\bf q}(s)^T,\nu(s))$.
It means that the definition of this kind of curves by the differential equation (\ref{eq6})
is automatically parametrized by the arclength \cite{jjo}.
Note that the point  $({\bf q}(s)^T,\nu(s))$ can lie outside the special characteristic surface of Section \ref{Kap8}.
But if we use the solutions of the eqs.(\ref{eq7}),(\ref{eq8}), the formulae (\ref{eq10}) and (\ref{eq11}), we get
again $D=0$ on the special characteristic surface $V({\bf q})=\nu$
\begin{equation}
 \int_o^s \left[ \sum_{k}^N (\frac{\partial V}{\partial q^k})^2(\tau) -G(\tau) \right]^{1/2} d\tau = 0 \ .
\label{eq17-2}
\end{equation}
$ \hfill{\rule{3mm}{3mm}}$ 
 
In the next subsection, we derivate the geodesic distance along a partial direction. 
This will usually be a direction which does not point along the steepest ascent line, or another geodesics.
We have for $j=1,...,N$, while $D \ne 0$,
\begin{equation} 
\frac{\partial}{\partial q^j} {\cal J}({\bf q},\nu ,s)= 
\frac{1}{2D^{1/2}}  \int_o^s \left[ 2 \sum_{k}^N \dot{q}^k(\tau) \frac{\partial \dot{q}^k}{\partial q^j}
  -( \frac{\partial}{\partial q^j} \frac{1}{G} ) \dot{\nu}^2   - \frac{2}{G} \dot{\nu}  \frac{\partial \dot{ \nu}}{\partial q^j}   \right] d\tau  \ .
\label{eq1p}
\end{equation}
With partial integration it becomes 
\begin{equation}
= 
\frac{1}{D^{1/2}} \left[  \sum_{k}^N \dot{q}^k(\tau) \delta^k_j   - \frac{1}{G} \dot{\nu} \frac{\partial \nu}{\partial q^j}  \right]_o^s
-\frac{1}{D^{1/2}} \int_o^s \left[ \frac{1}{2} (\frac{\partial}{\partial q^j} \frac{1}{G} ) \dot{\nu}^2
+\sum_{k}^N \ddot{q}^k(\tau) \frac{\partial q^k}{\partial q^j} + \ddot{\nu} \frac{\partial \nu}{\partial q^j}   \right] d\tau \ .
\label{eq21}
\end{equation}
It is equal
\begin{equation}
= 
\frac{1}{D^{1/2}} (\dot{q}^j(s)- \dot{q}^j(0)) 
-\frac{1}{D^{1/2}} \int_o^s \left[ \frac{1}{2} (\frac{\partial}{\partial q^j} \frac{1}{G} ) \dot{\nu}^2
+ \ddot{q}^j(\tau) \right] d\tau \ .
\label{eq22}
\end{equation}
The integrand of the right hand side summand is zero because of the geodesic property, Eq.(\ref{eq8}).
We get with (\ref{eq15})
\begin{equation}
\frac{\partial {\cal J}}{\partial q^j} =  \frac{\dot{q}^j}{{D^{1/2}}} =    \frac{\partial \sqrt{D}}{\partial \dot{q}^j} \ .
\label{eq4p}
\end{equation}
Analogously, we get
\begin{equation}
\frac{\partial {\cal J}}{\partial \nu} =  -\frac{\dot{\nu }}{{D^{1/2}\ G}} =   \frac{\partial \sqrt{D}}{\partial \dot{\nu}} \ .
\label{eq5p}
\end{equation}
We use Eq.(\ref{eq15}) for a further derivation \cite{cour}, p.119. 
Put $F=\sqrt{D}$. 
Then it is
\begin{equation}
 \sum_k^N \dot{q}^k \frac{\partial F}{\partial \dot{q}^k} + \dot{\nu} \frac{\partial F}{\partial \dot{\nu}}  = \frac{1}{F} (
 \sum_k^N (\dot{q}^k)^2  - \frac{1}{G} (\dot{\nu})^2 ) = F  \ .
\label{eq18}
\end{equation}
Now we have for $k=1,...,N$ using eqs.(\ref{eq21}) and (\ref{eq22})
\begin{equation}
{\cal J}_{q^k} = F_{\dot{q}^k} = \frac{1}{F} \sum_l^N \delta_k^l \dot{q}^k = \frac{\dot{q}^k}{F} \ , 
{\rm \ \ \  and \ \ \ } {\cal J}_{\nu } = F_{\dot{\nu} } = \frac{1}{F} \frac{(-1)}{G} \dot{\nu } \ . 
\label{eq19}
\end{equation}
The multiplication of eqs.(\ref{eq19}) with the inverse metric matrix $g^{JK}$ and 
a further multiplication with ${\cal J}_{q^j}  =   F_{\dot{q}^j}$ for $j=1,...,N$, and with ${\cal J}_{\nu } = F_{\dot{\nu}}$ 
and addition over all terms for  $j,k=1,...,N$, and $\nu$, results in
\begin{equation}
\sum_{k,j}^N  g^{kj} {\cal J}_{q^k} {\cal J}_{q^j} -G {\cal J}_{\nu }^2 = 
\sum_{k}^N {\cal J}_{q^k}^2 -G \, {\cal J}_{\nu }^2 =
\frac{1}{F} ( \sum_{j}^N F_{\dot{q}^j} \dot{q}^j + F_{\dot{\nu}} \dot{\nu } ) = 1  \ 
\end{equation}
because of eqs.(\ref{eq18}) and (\ref{eq19}). \ Now, put $\Lambda =  {\cal J}^2$ being the square of the geodesic distance 
between the points ${\bf q}$ and ${\bf q}_o$. 
In physics the ${\Lambda}$ is sometimes called the world function (with interchanged signature, and $\nu$ is replaced by the time).
With
\begin{equation}
\frac{\partial \Lambda}{\partial {q}^I} = 2 {\cal J} \ \frac{\partial{\cal J} }{\partial {q}^I}
\end{equation}
we have the well-known differential equation of the first order
\begin{equation}
\sum_{k,j}^N  g^{jk} \ {\Lambda}_{q^k} {\Lambda}_{q^j} -G \ ({\Lambda}_{\nu })^2 = \sum_{j}^N \ ({\Lambda}_{q^j})^2 -G \ ({\Lambda}_{\nu })^2 = 4\ \Lambda \  .
\label{eq44}
\end{equation}
\section{Singular part of the elementary solution}
\label{Kap10}
%
We use formula (2.10) of Czapor/McLenaghan \cite{mc} for the elementary solution of eq.\,(\ref{invari})
\begin{equation}
V({\bf q},\nu,{\bf q}_o,\nu_o) = \frac{1}{2 \pi } (\rho({\bf q},\nu,{\bf q}_o,\nu_o))^{-1/2}\ exp
\left( \frac{1}{4} \int_0^{s({\bf q},\nu)} \ A^I \ \Lambda_{,I} \frac{dt}{t} \right) 
\label{cza}
\end{equation}
where
\begin{equation}
 \rho({\bf q},\nu,{\bf q}_o,\nu_o)= 8 (\gamma({\bf q},\nu) \, \gamma({\bf q}_o,\nu_0) )^{1/2} \ 
 \left[ det(\frac{\partial^2 \Lambda }{\partial q^I \partial q_o^J } ) \right]^{-1} \ .
\end{equation}
The factor 8 here emerges because the Ref.\,\cite{mc} uses the dimension n=4 only. In general, we have the factor 
$2^{n-1}$. The part $\rho $ is the so called discriminant function, and $\gamma$ is the norm of the determinant of the 
metric (\ref{metric}).
One should also compare the theory in Ref.\,\cite{fried}, p.129 ff.
At the beginning we approximate the square of the geodesic distance $\Lambda$, see (\ref{eq17}).
We locally define the square of the geodisic distance by
\begin{equation} 
{\Lambda}({\bf q},\nu ;{\bf q}_o,\nu_o) =  \sum_{k}^N (q^k-q_o^k)^2  -\frac{1}{G({\bf q}_o)} (\nu -\nu_o)^2  + 
{\cal O}(| ({\bf q},\nu)- ({\bf q}_o,\nu_o)|^3)\ .
\label{eq75}
\end{equation}
It is ${\Lambda}({\bf q},\nu ;{\bf q}_o,\nu_o)=0$  the equation of the characteristic conoid with apex in $({\bf q}_o,\nu_o)$.
We can differentiate it. Up to linear terms it becomes 
\begin{equation}
\frac{\partial^2 \Lambda}{\partial q^i\,\partial q_o^j} = 
\left\{ \begin{array}{r@{\quad ;\quad}l}
          0 & 1 \le i,j \le N, i\ne j 
 \\      -2 & i=j  
        \end{array} 
\right. 
\end{equation}\
and
\begin{equation}
\frac{\partial^2 \Lambda}{\partial q^i\,\partial \nu_o} = 0, \quad
\frac{\partial^2 \Lambda}{\partial \nu\,\partial \nu_o^j} = \frac{2}{G({\bf q}_o)},  \quad
\frac{\partial^2 \Lambda}{\partial q_o^i\,\partial \nu} = \frac{4}{G({\bf q}_o)^2} (\nu -\nu_o){\cal H}_i , 
\end{equation}
The determinat of all entries is (if one weights all by factor 1/2)
\begin{equation}
det\left( \frac{1}{2}{\Lambda}_{q^I\,q_o^J} \right) = \frac{(-1)^N}{{G({\bf q}_o)}}  \  .
\end{equation}
Then the factor in formula (\ref{cza}) becomes
\begin{equation}
\frac{\left| det \left( \displaystyle\frac{1}{2} \displaystyle\frac{\partial^2 \Lambda}{\partial q^I\,\partial q_o^J} \right)\right|^{1/2} }
     {\gamma({\bf q})^{1/4} \gamma({\bf q}_o)^{1/4}} 
=  \left( \frac{G({\bf q})}{G({\bf q}_o)} \right)^{1/4} \ ,
\label{eq79}
\end{equation}
where we use the norm of the determinant, $\gamma=1/G$, of the metric (\ref{metric}).
The factor is really equal one at the beginning of the geodesic line, if $({\bf q},\nu)=({\bf q}_o,\nu_o)$. 

In the next step, we treat the integral of (\ref{cza}).
The factor (\ref{eq79}) only concerns the second order operator; but the integral also includes the first order terms.
We have to calculate 
\begin{equation}
exp\left( \frac{1}{4} \int_0^s ( \sum_i A^i\, \partial_i \Lambda + Tr {\bf H}\ \partial_{\nu} \Lambda )\frac{dt}{t} \right) \ 
\label{intg}
\end{equation}
where the path integral goes along a geodesic line from $({\bf q}_o,\nu_o)$ to $({\bf q},\nu)$.
We again use the explicite definition $\Lambda = {\cal J}^2$ and $ {\cal J}= F s$, see (\ref{eq17}) and (\ref{eq18}). The integrant is
\begin{equation}
 \sum_i \frac{{\cal H}_i}{G}\, 2 {\cal J} \frac{\dot{q}^i}{F} +  
Tr {\bf H}\, 2 {\cal J} \frac{(-\dot{\nu})}{F\,G} \ ,
\end{equation}
see (\ref{eq19}).  The integral in (\ref{intg}) becomes
\begin{equation}
\int_0^s \left( 2\,\sum_i \frac{{\cal H}_i}{G} \dot{q}^i - 2 Tr{\bf H} \frac{\dot\nu}{G} \right) dt \ 
\end{equation}
\begin{equation}
= \int_0^s  \left( - 2\sqrt{G} \sum_i  \frac{\partial}{\partial q^i}(\sqrt{\frac{1}{G}})\, \dot{q}^i - 2 Tr{\bf H}\frac{\dot\nu}{G} \right) dt \ .
\end{equation} 
The first sum is a gradient in the space of the {\bf q} variables,
but does not depend on $\nu$, we can write it as a pure derivation to $t$. 
Thus we get with the gradient theorem for path integrals 
\begin{equation}
=  \left. -2 log\,(\sqrt{\frac{1}{G}})\,\right|_{q_o}^q -2 \int_0^s Tr{\bf H} dt
=  log\,(\frac{G({\bf q})}{G({\bf q}_0)} )  -2 \int_0^s Tr{\bf H}\frac{d\nu}{G}  \ .
\end{equation}
The first term with factor 1/4 in (\ref{intg}) results in the same expression like 
the factor from the second order part in result (\ref{eq79}), 
and the last term remains a path integral giving the singular part of the elementary solution
\begin{equation}
U({\bf q},\nu;{\bf q}_0,\nu_0) =  \frac{1}{\sqrt{2} \pi} \left( \frac{G({\bf q})}{G({\bf q}_o)} \right)^{1/2}\  exp(\,-\frac{1}{2} \int_0^s \frac{Tr{\bf H}}{G} d\nu \,) \ . 
\label{ele}
\end{equation}
The elementary solution (or Green's function) is
\begin{equation}
U({\bf q},\nu;{\bf q}_0,\nu_0)\delta^{(\pm )}(\Lambda)+ W({\bf q},\nu;{\bf q}_0,\nu_0) \Theta^{(\pm )}(-\Lambda) \ .
\label{eleme}
\end{equation}
$\delta^{(\pm )}$ is the ($N+1$)-dimensional Dirac delta distribution, but
$\Theta^{(\pm )}$ is the Heaviside step function with support in the forward (+) or backward (-) conoid.
The functions $U$ and $W$ are unique. The equations for $W$ are hard to solve, and no analytical expression is known in a general metric.

Though the integrant of the path integral is independent on $\nu$ and $\nu_0$, 
the integral depends on the ($N$+first) variable, $\nu$, over the kind of pathway which determines the corresponding geodesic line. 
For example, two points $({\bf q}^T,\nu)$ and $({\bf q}_0^T,\nu_0)$ on a steepest ascent path are connected by the corresponding steepest ascent line,
see Section \ref{Kap7}.
However, if we compare the two points $({\bf q}^T,\nu_0)$ and $({\bf q}_0^T,\nu_0)$ in a pure ``space-like'' hyperplane,
$\nu_o$=constant, then we can choose $d\nu/ds=0$ and get with Eq.(\ref{eq8}) that the geodesic line is 
a straight line in the {\bf q}-space.
Thus, two space points  ${\bf q}$ and ${\bf q}_0$ are connected by two different geodesics depending on the value of the $\nu$ variable.

For the validity of a strong HP, it is to test  $P\,U \stackrel{?}{=} 0$ where $P$ is the adjoint differential operator to (\ref{waveop}). 
That would be immediately a sufficient HP property by $W=0$, see \cite{mc} as well as Ref.\,\cite{hadam}.  
Then only the first term in (\ref{eleme}) appears. It has support along the null conoid, and it represents a sharp, 
wake-free wave front.
The $W$-term is sometimes named the tail. It accounts for the possibility of a diffusive propagation. \ 
The adjoint differential operator $P$ to (\ref{waveop}) is, compare (\ref{invari})
\begin{equation}
P\,U= g^{I\,J} \, \nabla_I \, \nabla_J \ U({\bf q},\nu) - \, \nabla_I (A^I \ U({\bf q},\nu))  \ .
\label{PU}
\end{equation}
However, a calculation results in $P\,U \ne 0$, in general, as it is to expect from our conclusion 
in Section  \ref{Kap6}.
%
%
\section{Discussion}
\label{Kap11}
The present investigation can be regarded as an investigation of
operators that allow progressing wave families as defined by Courant. \cite{cour} 
The wave front (\ref{sol}) solves the differential equation (\ref{waveeq}) in an  $I\!\!R^{N+1}$. 
It satisfies a Huygens' principle if one uses a distribution ansatz in Eq.(\ref{sol}). 
By a Huygens' principle of the propagation of waves we first understand that sharp signals propagate as sharp signals.
The ''signal'' here is the $(N-1)$-contour to a corresponding level of the given surface V in  $I\!\!R^{N+1}$.
The differential equation (\ref{waveeq}) is linear, but its second order part has one nonconstant coefficient. 
Note that the nonlinear coefficients of operator (\ref{waveop}), G({\bf q}) and Tr {\bf H}({\bf q}), 
are built from the gradient and the Hessian of the function of interest, V({\bf q}). 
The operator (\ref{waveop}) is constructed for this function.
The function V({\bf q}) then also emerges in the very simple solution (\ref{sol}) being a progressing wave,  
for which we can also pose an appropriated Cauchy initial value problem (\ref{ivp_0}) and (\ref{ivp_1}). 
However, if we pose a general IVP for (\ref{waveeq}), 
we cannot expect such a simple solution like (\ref{sol}); 
in contrast, we will find a general solution formula, an integral over a full conoid region. 
Thus, the strong HP in the sense of Hadamard's 'minor premise' will not be fulfilled. 
That is the result of the calculation of diverse necessary criteria for the strong HP: they are not fulfilled.
We found that Eq.\,(\ref{waveeq}) itself is not a Huygens' equation. 
Wave operators that themselves satisfy the strong HP are known to be very rare, while those with families of 
progressing wave solutions are much more common. There may be no single connection between the two classes of equations \cite{bom91}.
Nevertheless, it is perhaps worth pointing out that operator $L$ in Eq.(\ref{waveeq}) allows 
undistorted progressing waves in every dimension,
and its splitting property (\ref{wavei}) is a hint that we could bring forward a Huygens' principle 
of the propagation of waves from
one kind of subdimension $M<N$, say odd, to the other kind, say even, and vice versa. 
Hadamards 'major premise' holds in any case, and for the special IVP (\ref{ivp_0}) and (\ref{ivp_1}) also 
the progressing wave solutions (\ref{sol}) are wake-free.
\\

If one compares operator\,(\ref{waveop}) with a wave equation with constant coefficients, then $1/G$ would
correspond to the velocity $c$ of the wave movement. In a certain sence this fits here as well:
If $G$ is small then the level hypersurfaces change by large steps, but if $G$ is large then the levels are dense. 
But note that this 'velocity' does not explicitely appear in the solution (\ref{sol}).
\\

Throughout the paper we have assumed that $G\ne 0$. The contrary case concerns the stationary points of the surface V of interest. 
There the 'rays' of operator (\ref{waveop}), the steepest ascent/steepest descent lines, meet in a singularity, 
and the level hypersurface reduces to a point (in minima or maxima), 
or it forms a caustic (in saddle points). 
The treatment of such saddle points is a further problem, see for example Ref.\,\cite{bena}. 
%
\section*{Acknowledgements}
Financial support from the Spanish Ministerio de Ciencia e Innovaci\'{o}n, DGI project CTQ2011-22505 and, 
in part from the Generalitat de Catalunya projects 2009SGR-1472 is fully acknowledged.
It is a pleasure to thank R.\,Schimming for a helpful discussion to Ref.\,\cite{boqu12}, 
and to thank M.\,Belger for hints to this paper.
Some formulas are controlled by a calculation with a Mathematica program of R.\,Sulanke \cite{sul13}.
%
%

%
\section*{Glossary}
Definitions (by Wikipedia) of field-specific terms in this paper:
\begin{itemize}
\item Characteristic surface $S$ for operator $L \psi =0$: if data for $\psi$ are prescribed on the surface $S$, 
then it may be possible to determine the normal derivative of  $\psi$ on $S$ from the differential equation. 
If the data on $S$ and the differential equation determine the normal derivative of  $\psi$ on $S$, then $S$ is non-characteristic. 
If the data on $S$ and the differential equation do not determine the normal derivative of  $\psi$ on $S$, 
then the surface is characteristic, and the differential equation restricts the data on $S$: 
the differential equation is internal to $S$.
\item Christoffel symbols: named for Elwin Bruno Christoffel (1829-1900), numerical arrays of real numbers that describe, in coordinates, 
the effects of parallel transport in curved surfaces and, more generally, manifolds.
\item Curvature tensor:  the infinitesimal geometry of Riemannian manifolds with dimension at least 3 
is too complicated to be described by a single number at a given point. 
Riemann introduced an abstract and rigorous way to define it, now known as the curvature tensor. 
\item Eikonal: the meaning of this equation is that any change in the region is pushed to move at right angles to the constant hypersurfaces, 
and along lines of steepest descent/ascent determined by the field of the vector and its sign.
\item Geodesic line: is a generalization of the notion of a straight line to curved spaces. In the presence of an affine connection, 
a geodesic is defined to be a curve whose tangent vectors remain parallel if they are transported along it. 
If this connection is the Levi-Civita connection induced by a Riemannian metric, 
then the geodesics are (locally) the shortest path between points in the space.
\item Huygen's principle: named for Christiaan Huygens (1629-1695), any point on a wave front may be regarded as the source of secondary waves. 
The surface that is tangent to the secondary waves can be used to determine the future position of the wave front.
\item Initial value problem:  (also called the Cauchy problem) is the differential equation together with specified values 
 of the unknown function and its derivative(s) at given points in the domain of the solution. 
\item Metric: named after Georg Friedrich Bernhard Riemann (1826-1866), the Riemannian metric $g$ is a symmetric tensor that is positive definite. 
It is an inner product on the tangent space at each point which varies smoothly from point to point. 
Riemannian geometry is the branch of differential geometry that studies smooth manifolds with a Riemannian metric.
This gives, in particular, local notions of angle, length of curves, surface area, and volume. 
From those some other global quantities can be derived by integrating local contributions.
\item Ricci tensor: named after Gregorio Ricci-Curbastro (1853-1925), represents the amount by which the volume of a geodesic ball in a curved Riemannian manifold 
deviates from that of the standard ball in Euclidean space. As such, 
it provides one way of measuring the degree to which the geometry determined by a given Riemannian metric might differ from that of ordinary Euclidean $I\!\!R^{N}$.
Like the metric itself, the Ricci tensor is a symmetric bilinear form on the tangent space of the manifold. 
\item Wave equation: a second-order linear partial differential equation for the description of waves – as they occur in physics – 
such as sound waves, light waves or water waves. Wave equations are examples of hyperbolic partial differential equations, 
but there are many variations.
\end{itemize}
\end{document}